\documentclass[a4paper, 10pt]{amsart}

\usepackage{amsmath,amsthm}
\usepackage{amssymb}

\usepackage{setspace}
\usepackage[latin1]{inputenc}
\usepackage{anysize} 

\marginsize{2cm}{2cm}{2cm}{3.5cm} 

\title{The extension of Buckley-Feuring solutions for non-polynomial fuzzy partial differential equations}

\author{D. Gálvez and J. L. Pino\\
{\tiny Departamento de Estad\'istica e Investigaci\'on
Operativa}\\
{\tiny Universidad de Sevilla}}

\date{\today}

\address{Facultad de Matemáticas. Universidad de
Sevilla, 41012 Sevilla, Spain}

\address{\textit{Keywords: Fuzzy differential equations, Buckley-Feuring solution, non-polynomial}}

\address{\textit{2000 Mathematics Subject Classification: 03E72,
46S40}}
\email{davidgalvez@us.es; jlpino@us.es}
\keywords{Fuzzy differential equations, Buckley-Feuring solution,
non-polynomial}

\subjclass{03E72, 46S40}

\begin{document}

 \maketitle

\vspace*{-0.8cm}

\begin{abstract}
This paper presents the natural extension of Buckley-Feuring
method proposed in \cite{BuckleyFeuring99} for solving fuzzy
partial differential equations (FPDE) in a non-polynomial
relation, such as the operator $\varphi(D_{x_1}, D_{x_2})$, which
maps to the quotient between both partials. The new assumptions
and conditions proceedings from this consideration are given in
this document.

\end{abstract}

\vspace*{-0.4cm}

\section*{Introduction}

Many approaches for obtaining non-numerical solutions of fuzzy
differential equations (FDE) have been developed from the
introduction of fuzzy set concept by Zadeh \cite{Zadeh65} . These
ones give a diversity of definitions for FDE solution based on
different notions of fuzzy derivative, such as Seikkala
derivative, Buckley-Feuring derivative, Puri-Ralescu derivative,
Kandel-Friedman-Ming derivative, Goetschel-Voxman derivative, or
Dubois-Prade derivative . Some relations between these derivatives
are presented by Buckley and Feuring in \cite{BuckleyFeuring00} .
However, only a few of these fuzzy derivatives are valid in some
contexts as FDE solution. For example, the Goetschel-Voxman
derivative, or the Dubois-Prade derivative provide solutions that
cannot be a fuzzy number, whereas Puri-Ralescu derivative, and
Kandel-Friedman-Ming derivative, always exist and provide a fuzzy
number as solution of the FDE, but making use of abstract
subtractions of fuzzy concepts in their definitions, making
difficult the interpretation of this solutions in some real
applications.

This paper uses the Buckley-Feuring derivative for solving FPDE.
This derivative does not always exist, but if it does, provides a
fuzzy number solution easily understandable in the context in
which a specific FPDE has been developed.

The authors which proposed this concept of derivative, developed a
methodology for solving constant coefficients polynomial FPDE in
\cite{BuckleyFeuring99} . This paper present the extension of this
methodology to a non-polynomial expression in partial fuzzy
derivatives.\\\\

In the following lines, the components of a FPDE are enumerated:
\begin{itemize}
\item $x_i,\quad i=1, 2,\quad
 x_1\in S_1\subset I_1=(0, M_1],\quad x_2\in S_2\subset I_2=(0, M_2]$. Other domain limits can be established in this subsets, such as $ x_1>x_2$.
\item $\tilde{\boldsymbol{\beta}}= (\tilde{\beta}_1,
\tilde{\beta}_2,...,\tilde{\beta}_k)$, a triangular fuzzy number
vector.
\item $\mu(\beta_j)$ is the membership function of $\beta_j\in \tilde{\beta_j}.$
\item $\mu_{\beta_j}(\alpha)=\{\beta_j\mid \mu(\beta_j)\geq\alpha,\quad\alpha\in(0,1)\}$
 set called $\alpha$-cut.\\\\ These sets are closed and bounded,
 so that is possible to define, for a fuzzy number
 $\tilde{\beta_j}: \tilde{\beta_j}[\alpha]=
 [b_1(\alpha), b_2(\alpha)]$, where:
 \begin{itemize}
 \item$b_1(\alpha)$ is the lower value $\beta_j$ in which
 $\mu(\beta_j)\geq \alpha,\quad \beta_j\in \tilde{\beta}_j$.

 \item$b_2(\alpha)$ is the higher value $\beta_j$ in which
 $\mu(\beta_j)\geq \alpha,\quad \beta_j\in \tilde{\beta}_j$.
 \end{itemize}

and $\tilde{\boldsymbol{\beta}}[\alpha]=
\prod_j\tilde{\beta_j}[\alpha]$

\item $\tilde{V}(x_1,x_2,\tilde{\boldsymbol{\beta}})$ is a positive and continuous function
in $(x_1, x_2)\in S_1 \times S_2$ with partials $D_{x_1},
D_{x_2}$. This function must be also strictly increasing or
strictly decreasing in $x_2\in S_2$, that is
$\tilde{V}(k,x_2,\tilde{\boldsymbol{\beta}})$ is strictly
increasing or strictly decreasing for all constant
$k\in\mathbb{R}$.\\The fuzzy character of
$\tilde{V}(x_1,x_2,\tilde{\boldsymbol{\beta}})$ shown by the tilde
placed over $V$ is fixed by $\tilde{\boldsymbol{\beta}}$, and
support the use of Buckley-Feuring derivative for solving FPDE.

\item $\varphi(D_{x_1}, D_{x_2})$ is an expression with constant coefficients in $(D_{x_1},
D_{x_2})$ applied  to
$\tilde{V}(x_1,x_2,\tilde{\boldsymbol{\beta}})$.

\item$F(x_1, x_2,\tilde{\boldsymbol{\beta}})$ continuous function in  $(x_1, x_2)\in
S_1 \times S_2.$

\end{itemize}

The specific FPDE treated in this paper has the following form
according with this notation:

$$\varphi(D_{x_1}, D_{x_2})\tilde{V}(x_1,x_2,\tilde{\boldsymbol{\beta}})= \frac{\partial \tilde{V}/\partial x_1}{\partial \tilde{V}/\partial
x_2}=F(x_1, x_2,\tilde{\boldsymbol{\beta}})$$

\section{The Buckley-Feuring (B-F) solution}
\label{sec:B-F}

The Buckley-Feuring (B-F) solution uses a solution of the crisp
partial differential equation:

$$V(x_1,x_2)= G(x_1,x_2,\boldsymbol{\beta}),$$
with $G$ continuous $\forall (x_1,x_2)\in S_1\times S_2$.\\

The next step is the fuzzification of $G$:
$$\tilde{Y}(x_1,x_2)= \tilde{G}(x_1,x_2,\tilde{\boldsymbol{\beta}}),$$
with $\tilde{G}$ continuous $\forall (x_1,x_2)\in S_1 \times S_2$
and strictly monotone for $x_2\in S_2$. Note that $\tilde{Y}_i$ is
only the fuzzy representation of $G$, but not necessary the
solution of the fuzzy partial differential equation. If it finally
happens and $\tilde{Y}(x_1,x_2)$ is a B-F solution, then
$\tilde{V}(x_1,x_2,\tilde{\boldsymbol{\beta}})=\tilde{Y}(x_1,x_2)$.

With this notation, it is possible to see that:\\\\
$\tilde{Y}(x_1,x_2)[\alpha]= [y_1(x_1,x_2,\alpha),
y_2(x_1,x_2,\alpha)]$, and\\\\
$\tilde{F}(x_1,x_2,\tilde{\boldsymbol{\beta}})[\alpha]=
[f_1(x_1,x_2,\alpha), f_2(x_1,x_2,\alpha)], \forall\alpha$.\\

and, by definition:

$y_1(x_1,x_2,\alpha)= \min\{G(x_1,x_2,\boldsymbol{\beta}),\quad
\boldsymbol{\beta}\in\tilde{\boldsymbol{\beta}}[\alpha]\}$,\\
$y_2(x_1,x_2,\alpha)= \max\{G(x_1,x_2,\boldsymbol{\beta}),\quad
\boldsymbol{\beta}\in\tilde{\boldsymbol{\beta}}[\alpha]\}$
,and\\\\
$f_1(x_1,x_2,\alpha)= \min\{F(x_1,x_2,\boldsymbol{\beta}),\quad
\boldsymbol{\beta}\in\tilde{\boldsymbol{\beta}}[\alpha]\}$,\\
$f_2(x_1,x_2,\alpha)= \max\{F(x_1,x_2,\boldsymbol{\beta}),\quad
\boldsymbol{\beta}\in\tilde{\boldsymbol{\beta}}[\alpha]\}$,\\
$\forall x_1, x_2, \alpha$.\\

If it is possible to apply the $\varphi(D_{x_1}, D_{x_2})$
operator to $y_i, i=1, 2 $, getting continuous expressions
$\forall(x_1, x_2)\in S_1\times S_2,\forall\alpha$, then it will
be feasible to define the following expression in this domain
$\Gamma(x_1,x_2,\alpha)$:

$$\Gamma(x_1,x_2,\alpha)=[\Gamma_1(x_1,x_2,\alpha),\Gamma_2(x_1,x_2,\alpha)]$$

with:

$\Gamma_1(x_1,x_2,\alpha)= \varphi(D_{x_1},
D_{x_2})y_1(x_1,x_2,\alpha)$

$\Gamma_2(x_1,x_2,\alpha)=\varphi(D_{x_1},
D_{x_2})y_1(x_1,x_2,\alpha)$\\\\

For a B-F solution $Y$, it is necessary to be a fuzzy number for
this one. If, for each pair $(x_1, x_2)\in S_1\times S_2$,
$\Gamma(x_1,x_2,\alpha)$ defines an $\alpha$-cut of a fuzzy
number, then, Buckley and Feuring \cite{BuckleyFeuring99} call to
$\tilde{Y_i}(x_1,x_2)$ \textbf{differentiable}, and we can write:
$$\varphi(D_{x_1},
D_{x_2})\tilde{Y}(x_1,x_2)[\alpha]=\Gamma(x_1,x_2,\alpha),$$
$\forall(x_1, x_2)\in S_1\times S_2,\quad\forall\alpha$.\\\\

So that, it is necessary to test if  $\Gamma(x_1,x_2,\alpha)$
really define an $\alpha$-cut of a fuzzy number and verify the
differentiability of $\tilde{Y}(x_1,x_2)$. For a triangular fuzzy
number, the conditions are \cite{GoetschelVoxman86} :

\begin{enumerate}
\item$\varphi(D_{x_1},
D_{x_2})y_1(x_1,x_2,\alpha)$ is an increasing function of
$\alpha$, for each $(x_1, x_2)\in S_1\times S_2$.
\item$\varphi(D_{x_1},
D_{x_2})y_2(x_1,x_2,\alpha)$ is a decreasing function of $\alpha$,
for each $(x_1, x_2)\in S_1\times S_2$.
\item$\varphi(D_{x_1},
D_{x_2})y_1(x_1,x_2,1)\leq\varphi(D_{x_1}, D_{x_2})y_2(x_1,x_2,1)$
for each $(x_1, x_2)\in S_1\times S_2$.
\end{enumerate}

Once delimited the differentiability concept of
$\tilde{Y}(x_1,x_2)$, it is possible to define the Buckley-Feuring
solution. $\tilde{Y_i}(x_1,x_2)$ is a Buckley-Feuring  solution if
the following conditions hold:\\

\begin{enumerate}
\item $\tilde{Y_i}(x_1,x_2)$ is differentiable.
\item $\varphi(D_{x_1},
D_{x_2})\tilde{Y}(x_1,x_2)= \tilde{F}(x_1,
x_2,\boldsymbol{\tilde{\beta}})$.\\

\end{enumerate}

Obviously, if differentiability conditions hold by the candidate
to B-F solution $\tilde{Y}(x_1,x_2)$, this one will be a fuzzy
number. To complete the conditions for B-F solutions only is
necesary to test that:

$$\varphi(D_{x_1},
D_{x_2})\tilde{Y}(x_1,x_2)= \tilde{F}(x_1,
x_2,\tilde{\boldsymbol{\beta}}),$$

or the equivalent condition:

\begin{enumerate}
\item $\varphi(D_{x_1},
D_{x_2})y_1(x_1,x_2,\alpha)= F_1(x_1, x_2,\alpha).$
\item $\varphi(D_{x_1},
D_{x_2})y_2(x_1,x_2,\alpha)= F_2(x_1, x_2,\alpha).$
\end{enumerate}
$\forall(x_1, x_2)\in S_1\times S_2,\quad\forall\alpha$.\\\\

In this case we can identify $\tilde{Y}(x_1,x_2)$ with
$\tilde{V}(x_1,x_2,\tilde{\boldsymbol{\beta}})$.

\section{Boundary conditions}
\label{sec:Bound. cond.}

The FPDE can be subject to certain boundary conditions in a big
variety of forms depending on a constant vector
$\boldsymbol{c}=(c_1,...c_n)\in C_1\times ...\times C_n$. The
inclusion of this ones has not great consequences in the
methodology exposed.

The crisp solution $G$ acquires under boundary conditions the form
$G(x_1,x_2,\boldsymbol{\beta},\boldsymbol{c})$. The fuzzification
of $G$ can take $\boldsymbol{c}$ in a triangular fuzzy vector
$\tilde{\boldsymbol{c}}=(\tilde{c_1},...\tilde{c_n})\in C_1\times
...\times C_n$ with $\tilde{\boldsymbol{c}}[\alpha]=
\prod_i\tilde{C_i}[\alpha],\quad i=1,...n$
 and:

$$\tilde{Y}(x_1,x_2)= \tilde{G}(x_1,x_2,\tilde{\boldsymbol{\beta}},\tilde{\boldsymbol{c}})$$

In this  environment with boundary conditions, it is necessary to
add a new condition for a B-F solution: $\tilde{Y}(x_1,x_2)$ must
satisfier these conditions. In this form $\tilde{Y_i}(x_1,x_2)$ is
a Buckley-Feuring  solution if the following conditionshold:

\begin{enumerate}
\item $\tilde{Y_i}(x_1,x_2)$ is differentiable.
\item $\varphi(D_{x_1},
D_{x_2})\tilde{Y}(x_1,x_2)= \tilde{F}(x_1,
x_2,\tilde{\boldsymbol{\beta}})$.
\item$\tilde{Y}(x_1,x_2)$
satisfies the boundary conditions.
\end{enumerate}

\section{Example}
\label{sec:example}

The following FPDE is proposed:

$$\dfrac{\partial {\tilde{V}}/\partial{x_1}}{\partial \tilde{V}/\partial
x_2}= \tilde{\beta} x_1^{ -1}x_2,\quad \tilde{\beta}\in (0,
1),x_1\geq 1, x_2>0
$$

In this example:

$\tilde{F}(x_1, x_2,\boldsymbol{\tilde{\beta}})=\tilde{\beta}
x_1^{ -1}x_2,\quad \tilde{\beta}\in (0, 1),x_1\geq 1, x_2>0$\\

and, by definition, the operator\\

 $\varphi(D_{x_1}, D_{x_2})V(x_1,x_2) \longrightarrow\dfrac{\partial \tilde{V}/\partial x_1}{\partial
\tilde{V}/\partial x_2},$\\\\

A possible solution to this FPDE in a crisp environment is,
without special boundary conditions:

$$G(x_1,x_2,\boldsymbol{\beta})= x_1^{\beta}x_2+\gamma,\beta\in (0, 1),x_1\geq 1, x_2>0$$

with $\boldsymbol{\beta}=(\beta,\gamma)$

Applying the fuzzification in $\beta$ and $\gamma$, these ones
acquire a triangular fuzzy number form and
$\boldsymbol{\tilde{\beta}}
=(\tilde{\beta},\tilde{\gamma}),\quad\tilde{\beta}\in (0, 1)$.
While $G$ holds all the conditions required, the Buckley-Feuring
solution candidate is:

$$\tilde{Y}(x_1,x_2)=G(x_1,x_2,\tilde{\beta},\tilde{\gamma})= x_1^{\tilde{\beta}}x_2+\tilde{\gamma},\tilde{\beta}\in (0, 1),x_1\geq 1, x_2>0$$

The fuzzy parameters have a membership function associated
$\mu(\beta)$ and $\mu(\gamma)$ respectiveness. From the
$\alpha$-cuts, it is possible to define:

 $\tilde{\beta}[\alpha]=
 [b_1(\alpha), b_2(\alpha)]$,\\
 $\tilde{\gamma}[\alpha]=
 [g_1(\alpha), g_2(\alpha)]$, and\\
 $\tilde{\boldsymbol{\beta}}[\alpha]=
 \tilde{\beta}[\alpha] \times \tilde{\gamma}[\alpha]$\\

And from this ones:

$\tilde{Y}(x_1,x_2)[\alpha]= [y_1(x_1,x_2,\alpha),
y_2(x_1,x_2,\alpha)]$, and\\
$\tilde{F}(x_1,x_2,\tilde{\boldsymbol{\beta}})[\alpha]=
[f_1(x_1,x_2,\alpha), f_2(x_1,x_2,\alpha)], \forall\alpha $.\\

Where:

$y_1(x_1,x_2,\alpha)= \min\{G(x_1,x_2,\boldsymbol{\beta}),\quad
\boldsymbol{\beta}\in\boldsymbol{\tilde{\beta}}[\alpha]\}= G(x_1,x_2,g_1(\alpha), b_1(\alpha)) $,\\
$y_2(x_1,x_2,\alpha)= \max\{G(x_1,x_2,\boldsymbol{\beta}),\quad
\boldsymbol{\beta}\in\boldsymbol{\tilde{\beta}}[\alpha]\}=G(x_1,x_2, g_2(\alpha), b_2(\alpha))$ \\\\and,\\\\
$f_1(x_1,x_2,\alpha)= \min\{F(x_1,x_2,\boldsymbol{\beta}),\quad
\boldsymbol{\beta}\in\boldsymbol{\tilde{\beta}}[\alpha]\}=F(x_1,x_2, g_1(\alpha), b_1(\alpha))$,\\
$f_2(x_1,x_2,\alpha)= \max\{F(x_1,x_2,\boldsymbol{\beta}),\quad
\boldsymbol{\beta}\in\boldsymbol{\tilde{\beta}}[\alpha]\}= F(x_1,x_2, g_2(\alpha), b_2(\alpha))$,\\\\
$\forall\alpha\quad x_1\geq
1, x_2>0  $\\\\

In the $\tilde{G}$ function proposed, $\forall(x_1, x_2)x_1\geq
1, x_2>0$:\\

$y_1(x_1,x_2,\alpha)=x_1^{b_1(\alpha)}x_2+{g_1(\alpha)},\quad b_1(\alpha)\in (0, 1)$\\
$y_2(x_1,x_2,\alpha)=x_1^{b_2(\alpha)}x_2+{g_2(\alpha)},\quad b_2(\alpha)\in (0, 1)$\\\\
$f_1(x_1,x_2,\alpha)=b_1(\alpha)x_1^{-1}x_2,\quad b_1(\alpha)\in
(0,
1)$\\
$f_2(x_1,x_2,\alpha)=b_2(\alpha)x_1^{-1}x_2,\quad b_2(\alpha)\in
(0,
1)$\\\\

Testing the differentiability of
$\tilde{G}(x_1,x_2,\tilde{\boldsymbol{\beta}})$, from
$\Gamma(x_1,x_2,\alpha)=[\varphi(D_{x_1},
D_{x_2})y_1(x_1,x_2,\alpha), \varphi(D_{x_1},
D_{x_2})y_2(x_1,x_2,\alpha)]$ verifying if
$\Gamma(x_1,x_2,\alpha)$ defines an $\alpha$-cut of a triangular
fuzzy number for each pair $(x_1, x_2), x_1\geq 1, x_2>0$, the
following conditions must be hold:

\begin{enumerate}
\item$\varphi(D_{x_1},
D_{x_2})y_1(x_1,x_2,\alpha)= b_1(\alpha)x_1^{-1}x_2,\quad
b_1(\alpha)\in (0, 1)$ is an increasing function of  $\alpha$, for
each pair $(x_1, x_2), x_1\geq 1, x_2>0$.

It will happen if $\varphi(D_{x_1}, D_{x_2})y_1(x_1,x_2,\alpha)$
has positive derivative on $\alpha$:

$$d_{\alpha}\left(b_1(\alpha)x_1^{-1}x_2\right)=
b_1(\alpha)'x_1^{-1}x_2$$

While $\beta$ is a triangular fuzzy number and $b_1(\alpha)$is
defined from its $\alpha-cuts$, $b_1(\alpha)$ satisfies the
condition and is an increasing function with $b_1(\alpha)'>0$, and
$x_1^{-1}x_2>0$ for $ x_1\geq 1, x_2>0$,
$d_{\alpha}\left(b_1(\alpha)x_1^{-1}x_2\right)>0$ and the
condition holds.\\\\

\item$\varphi(D_{x_1},
D_{x_2})y_2(x_1,x_2,\alpha)= b_2(\alpha)x_1^{-1}x_2,\quad
b_2(\alpha)\in (0, 1)$ is a decreasing function of  $\alpha$, for
each pair $(x_1, x_2), x_1\geq 1, x_2>0$.

Again, it will happen if $\varphi(D_{x_1},
D_{x_2})y_2(x_1,x_2,\alpha)$ has negative derivative on $\alpha$:

$$d_{\alpha}\left(b_2(\alpha)x_1^{-1}x_2\right)=
b_2(\alpha)'x_1^{-1}x_2$$

We can now use that $\beta$ is a triangular fuzzy number and
$b_2(\alpha)$ is defined from its $\alpha-cuts$, so that, in
analogy, $b_2(\alpha)$ satisfies the condition and is a decreasing
function with $b_2(\alpha)'<0$, and $x_1^{-1}x_2>0$ for $ x_1\geq
1, x_2>0$, $d_{\alpha}\left(b_1(\alpha)x_1^{-1}x_2\right)<0$ and
the condition holds.\\\\

\item$\varphi(D_{x_1},
D_{x_2})y_1(x_1,x_2,1)\leq \varphi(D_{x_1}, D_{x_2})y_2(x_1,x_2,1)
\forall x_1, x_2)\in I_1\times I_2$.\\In the example:
$$b_1(1)x_1^{-1}x_2
\leq
b_2(2)x_1^{-1}x_2$$\\

This condition is holds automatically if we realize that
$b_1(\alpha)$ and $b_2(\alpha)$ are defined from a triangular
fuzzy number and $b_2(1)>b_1(1)$

\end{enumerate}

At this point, it is possible to say that
$\tilde{Y}(x_1,x_2)=G(x_1,x_2,\tilde{\beta},\tilde{\gamma})=
x_1^{\tilde{\beta}}x_2+\tilde{\gamma},\tilde{\beta}\in (0,
1),x\geq 1$ is differentiable and good candidate for B-F solution.
But it is necessary the last step, and it must satisfier that:\\\\

$\varphi(D_{x_1}, D_{x_2})\tilde{Y}(x_1,x_2)= \tilde{F}(x_1,
x_2,\tilde{\boldsymbol{\beta}}) \quad \beta, \gamma\in (0, 1),
\forall(x_1, x_2), x_1\geq 1, x_2>0$.\\

For $G(x_1,x_2,\tilde{\beta},\tilde{\gamma})=
x_1^{\tilde{\beta}}x_2+\tilde{\gamma},\tilde{\beta}\in (0,
1),x\geq
1, x_2\geq 0$:\\\\
$\varphi(D_{x_1}, D_{x_2})\tilde{Y}(x_1,x_2)= \tilde{\beta} x_1^{
-1}x_2,\quad \tilde{\beta}\in (0, 1),x_1\geq 1, x_2>0 $\\

and,\\

$\tilde{F}(x_1, x_2,\tilde{\boldsymbol{\beta}})=\tilde{\beta}
x_1^{ -1}x_2,\quad \tilde{\beta}\in (0, 1),x_1\geq 1, x_2>0.
$\\

Thus $G(x_1,x_2,\tilde{\beta},\tilde{\gamma})=
x_1^{\tilde{\beta}}x_2+\tilde{\gamma},\tilde{\beta}\in (0,
1),x\geq 1, x_2\geq 0$ is a Buckley-Feuring solution for this
non-polynomial fuzzy partial differential equation.


\begin{thebibliography}{99}


\bibitem{Allahviranloo02} Allahviranloo, T. Diference methods for fuzzy partial differential equations, Computational Methods in Applied Mathematics, Vol.2 (nº3) (2002), pp.
233-242.


\bibitem{BuckleyFeuring99} Buckley, J.J. and Feuring,T. Introduction to
fuzzy partial differential equations, Fuzzy Sets and Systems, 105
(1999), pp.241-248.


\bibitem{BuckleyFeuring00}  Buckley, J.J. and Feuring,T. Fuzzy  differential equations, Fuzzy Sets and Systems,
110 (2000), pp.43-54.

\bibitem{GoetschelVoxman86} Goetschel, R. and Voxman, W. Elementary fuzzy calculus, Fuzzy Sets and Systems,
18 (1986), pp.319-330.

\bibitem{Quan92} Quan, H., Hua, Y. and Jones, J.D. A general
method for calculating functions of fuzzy numbers, Applied
Mathematics Letters, Vol.5 (nº6) (1992), pp.51-55.



\bibitem{Seikkala87} Seikkala, s. On the fuzzy initial value problem, Fuzzy Sets and Systems,
24 (1987), pp.31-43.

\bibitem{Zadeh65} Zadeh, L.A. Fuzzy Sets. Inf. Control, 8 (1965),
pp.338-353.


\end{thebibliography}
\end{document}